\documentclass{amsart}
\usepackage[utf8]{inputenc}

\usepackage{graphicx}        
\usepackage{multicol}        
\usepackage[bottom]{footmisc}
\usepackage{amsfonts, amsmath, amssymb}
\usepackage{booktabs} 
\usepackage{enumitem}
\usepackage{comment}
\usepackage[linesnumbered,ruled,vlined]{algorithm2e}
\usepackage[draft=true]{minted} 
\setminted{breaklines=true, breakanywhere=true}
\setminted{fontsize=\footnotesize}

\usepackage{verbatim}
\usepackage[backend=bibtex,style=alphabetic,doi=true,isbn=false,url=false,giveninits=true,maxbibnames=50]{biblatex}
\usepackage{enumitem}
\setenumerate{label=(\roman*)}

\usepackage{placeins}
\usepackage{hyperref}

\newcommand{\nc}{\newcommand}
\nc{\lie}[1]{\mathfrak{#1}}
\nc{\g}{\mathfrak{g}}
\nc{\ra}{\rightarrow}
\nc{\ssl}{\mathfrak{sl}}
\nc{\ve}{\varepsilon}
\nc{\vp}{\varphi}
\nc{\bl}{\lambda}
\nc{\C}{\mathbb{C}}
\nc{\CC}{\mathbb{C}}
\nc{\OSCAR}{\text{OSCAR }}
\nc{\GAP}{\text{GAP }}
\nc{\polymake}{\text{polymake }}
\nc{\julia}{\text{julia }}

\newcommand{\ba}{\mathbf{a}}
\newcommand{\bb}{\mathbf{b}}

\newcommand{\bm}{\mathbf{m}}

\newcommand{\bv}{\mathbf{v}}

\theoremstyle{plain}
\newtheorem{theorem}[subsection]{Theorem}
\newtheorem*{theorem*}{Theorem}
\theoremstyle{definition}

\newtheorem{conjecture}{Conjecture}

\newtheorem{example}[subsection]{Example}
\newtheorem{definition}[subsection]{Definition}

\title{Computing monomial bases in Lie theory using OSCAR }
\author[]{Xin Fang and Ghislain Fourier and Lars Göttgens and Ben Wilop}
\address{Chair of Algebra and Representation Theory, RWTH Aachen University}
\email{fang@art.rwth-aachen.de}
\email{fourier@art.rwth-aachen.de}
\email{goettgens@art.rwth-aachen.de}
\email{ben@wilop.de}

\begin{document}

\begin{abstract}
In this survey, we present a detailed guide on using the computer algebra system \OSCAR to compute monomial bases for simple, finite-dimensional modules of simple, complex Lie algebras. We will also demonstrate how to determine monomial bases for the homogeneous coordinate ring of a (partial) flag variety, depending on a chosen birational sequence and a monomial order. This survey will be updated to reflect any advancements in OSCAR's capabilities in these areas.
\end{abstract}
\maketitle

\section*{Introduction}
For a finite-dimensional complex simple Lie algebra $\g$ and its finite-dimensional irreducible representation $V(\lambda)$,
finding an explicit vector space basis of $V(\lambda)$ is one of the most important questions in the representation theory.
The study of this question dates back to the work of Young, Hodge (Young tableaux) in the computation of the Hilbert polynomial of a Grassmann variety in its Plücker embedding;
 and the work of Gelfand-Tsetlin in the construction of such a basis for type $A$ Lie algebras with explicit formulae for the action of Chevalley generators. Recent developments around this question include the work of Lakshmibai-Seshadri-Littelmann, Chirir\`i-Fang-Littelmann \cite{LS5,LS98,CFL22,CFL23} (standard monomial basis), Lusztig \cite{Lus91} (canonical basis, semi canonical basis), Kashiwara \cite{Kas90} (global crystal basis), Mirkovi\'c-Vilonen \cite{MV} (MV basis), Feigin-Fourier-Littelmann \cite{FFL1,FFL2} (FFLV basis), Gross-Hacking-Keel-Kontsevich \cite{GHKK} (theta basis), to name but a few.

Many of such bases defined and studied in the past three decades come from different geometrizations of the universal enveloping algebra and their representations, making the computation of these bases usually a hard question. In \cite{FaFoL} the first two authors, together with Peter Littelmann, introduced a general construction of a family of monomial bases which have both theoretical and computational convenience. The computation of these bases has already been implemented in  \OSCAR\hspace{-1mm}.

\medskip
Birational sequences are introduced in \cite{FaFoL} as a first approximation of some of the above-mentioned bases.
From the Borel-Weil-Bott theorem, such bases (or rather their dual bases) can be considered as functions on the flag varieties or their Schubert subvarieties.
The key idea of the approach in \emph{loc.cit.} is to use a birational sequence to give parametrizations of the basis elements by polynomials.
By fixing a monomial order of the variables in the polynomials, we chose a leading monomial to represent the corresponding basis element.
This procedure is similar to replacing an algebra by its initial algebra.
Passing to the leading monomial is the first approximation of the basis element, and information will definitely get lost.
What one can still read off from the leading monomial is a skeleton of the flag variety or the Schubert variety.
In the mathematical language, this construction produces a flat family over $\mathbb{C}$ with generic fibre the flag variety or the Schubert variety,
and the special fibre is a reduced toric scheme.

The main conjecture is:
\begin{conjecture}\label{cnj:main}
    The special fibre is for any chosen birational sequence and monomial ordering an affine toric variety.
\end{conjecture}

The notion of birational sequences is motivated by the Newton-Okounkov theory \cite{KK,LM} and the PBW filtration \cite{FFL3}. In this language, the first approximation mentioned above is given by exactly the Newton-Okounkov body associated to the birational sequence and the monomial order, and the conjecture can be reformulated as: the associated Newton-Okounkov body is a rational polytope.

Birational sequences are used as a theoretical tool to unify many known toric degenerations of flag varieties and Schubert varieties. That is not all. It opens the door of computing Newton-Okounkov bodies for flag varieties and Schubert varieties, and the key idea is the notion of essential monomials, transforming the computation of the leading monomial into a down-to-earth and explicit representation theory problem. It is this notion which allows us to execute the computation in \OSCAR\hspace{-1mm}, and this survey serves as a guide for both theoretical and practical aspects.

The input data are ($N$ being the number of positive roots of $\lie g$)
\begin{itemize}
    \item a birational sequence $S$ (introduced in \cite{FaFoL}), for our purpose: a sequence of length $N$ of (not necessarily distinct) positive roots, such that ordered monomials form a generating set of the universal enveloping algebra of the positive part of $\lie g$,
    \item a fixed monomial ordering $<$ on $\mathbb{Z}^N$,
    \item a dominant, integral weight $\lambda$.
\end{itemize}
The \OSCAR package then computes a monomial basis of $V(\lambda)$ with respect to the generating set (we consider monomials in the root vectors of the corresponding negative roots of $S$) and the chosen monomial order. If we denote this set of essential monomials $\text{es}(S, <, \lambda)$, then $\text{es}(S, <, \lambda) + \text{es}(S, <, \mu) \subseteq \text{es}(S, <, \lambda + \mu)$ (where the sum is the Minkowski sum of the exponents). Consequently, \OSCAR computes then a minimal (with respect to a graded lexicographic ordering on weights) tuple $(\mu_1, \ldots, \mu_s)$ with $\mu_i \leq \lambda$ such that $\lambda = \mu_1 + \ldots+ \mu_s$ and $\text{es}(S, <, \lambda) = \sum a_i \text{es}(S, <, \mu_i)$ for some $a_i \geq 0$. The ``best cases'' are those, where one needs fundamental weights only and obtains any set of essential monomials through Minkowski sums.

We will extend the functionality soon to Demazure modules and Schubert varieties, and on the other hand connect to the realm of \OSCAR by using the capabilities provided by \polymake to actually compute the Newton-Okounkov bodies. The initial version of this survey has been published as \cite{OSCAR-book-chapter}, while this version will be extended by the time.

\vspace{10mm}
\noindent
\textbf{Organization of the survey}:
\begin{enumerate}[label=\arabic*.]
    \item Short introduction to \OSCAR\hspace{-1mm}.
    \item Theoretical background on birational sequences, essential monomials and Newton-Okounkov bodies.
    \item The implemented algorithm is explained.
    \item A lot of example codes for computing (important) monomial bases.
    \item How to compute bases of the homogeneous coordinate ring.
    \item Some experiments on generators of the monoid.
    \item Comparison of runtime against \GAP\hspace{-1mm}.
\end{enumerate}

\medskip

\noindent\textbf{Acknowledgements:} The work of GF, LG, and BW is funded by the Deutsche Forschungsgemeinschaft (DFG, German Research Foundation): ``Symbolic Tools in Mathematics and their Application'' (TRR 195, project-ID 286237555).

\section{\OSCAR}

\OSCAR \cite{OSCAR,OSCAR-book} stands for ``\textbf{O}pen \textbf{S}ource \textbf{C}omputer \textbf{A}lgebra \textbf{R}esearch'' and is a comprehensive computer algebra system for computations in algebra, geometry, and number theory.
It is built on top of multiple already existing computer algebra systems (GAP, Singular, polymake, Antic) that are specialized in one specific subject, and combines these underlying systems to connect these different topics for a user. 
The main focus of \OSCAR is to provide a consistent interface that delegates transparently to the underlying systems, while keeping the efficiency of already existing specialized algorithms. This is achieved by adding a semantical layer on top of the glue code, as well as functionality in more subjects.
Most parts of \OSCAR are written in the julia language \cite{Julia-2017}.

For the installation of \OSCAR\!\!, we refer the reader to the installation instructions\footnote{\url{https://www.oscar-system.org/install/}}. The features presented in this work are expected to be released with \OSCAR version 1.1. For access prior to that date, use any development version from the \texttt{master} branch of \OSCAR\!\!\footnote{\url{https://github.com/oscar-system/Oscar.jl}} as of \texttt{2024-03-21}\footnote{git commit hash \texttt{e4980cf945a96856def3650189a05166d7eb1502 }} or newer.
To install a non-released version of \OSCAR\!\!, replace step 3 in the installation instructions by
\begin{minted}{jlcon}
using Pkg
Pkg.add(url="https://github.com/oscar-system/Oscar.jl.git", rev="<git-ref>")
\end{minted}
\noindent where \texttt{<git-ref>} is either \texttt{master} or a git commit hash.

\section{Theoretical background}
The standard reference on basics of semisimple Lie algebras and their representation theory is \cite{Hum}.

Motivated by the theory of Newton-Okounkov bodies and the PBW filtration, studying monomial bases of representations serves as an intermediate step towards understanding monomial bases of finite-dimensional irreducible representations of finite-dimensional complex semisimple Lie algebras.
Geometrically it has a strong connection to the construction of toric degenerations of (partial) flag varieties.
Roughly speaking, the idea is to replace a complicated basis element, written in a sequence of generators of the Lie algebra $\g$, by a well-chosen leading term, and to express the leading term as a monomial of the fixed generators.

\subsection{Representations as functions}
Let $G$ be a simply connected semisimple algebraic group with Lie algebra $\g$.
Fixing a triangular decomposition $\g=\mathfrak{n}^+\oplus\mathfrak{h}\oplus\mathfrak{n}^-$ of $\g$ gives a Borel subgroup $B$ and a maximal torus $T$ of $G$ with Lie algebras $\mathfrak{n}^+\oplus \mathfrak{h}$ and $\mathfrak{h}$ respectively. The complete flag variety $G/B$ can be embedded into $\mathbb{P}(V(\lambda))$, for $\lambda$ a regular dominant weight, as a highest weight orbit, from which it gets the structure of a projective variety. The homogeneous coordinate ring of $G/B$ with respect to this embedding is then given by:
$$\mathbb{C}[G/B]\cong \bigoplus_{k\geq 0} V(k\lambda)^*,$$
where non-zero functions in $V(k\lambda)^*$ are homogeneous of degree $k$.

The existence of certain nice monomial bases for all $V(k\lambda)$ allows us to degenerate $\mathbb{C}[G/B]$ to its initial algebra, which is a toric algebra.
In \cite{FaFoL}, a general procedure is introduced to provide a systematical way of constructing monomial bases of $V(\lambda)$.
The key ingredient is the notion of a birational sequence, as will be explained in the next section.

\subsection{Birational sequences}

Let $U^-$ be the unipotent subgroup of $G$ with Lie algebra $\mathfrak{n}^-$ of dimension say $N$.
For each positive root $\beta$ of $\g$, we fix a non-zero root vector $f_\beta\in\mathfrak{n}^-$ of weight $-\beta$.
The corresponding root subgroup is denoted by $U_{-\beta}:=\{\exp(t f_{\beta})\mid t\in\mathbb{C}\}\subseteq U^-$.

\begin{definition}
A \emph{birational sequence} $(\beta_1,\beta_2,\ldots,\beta_N)$ is a sequence of positive roots (repetitions are allowed) of $\g$ such that the multiplication map
$$U_{-\beta_1}\times\cdots\times U_{-\beta_N}\dashrightarrow U^-,\ \ (u_1,\ldots,u_N)\mapsto u_1\cdots u_N$$
is a birational map.
\end{definition}

\begin{example}\label{Ex:Bir}
Let $\underline{w}_0=s_{i_1}s_{i_2}\cdots s_{i_N}$ be a reduced decomposition of $w_0$, the longest element in the Weyl group of $\g$.
\begin{enumerate}
\item[(1)] The sequence $(\alpha_{i_1},\alpha_{i_2},\ldots,\alpha_{i_N})$ is birational.
\item[(2)] For $1\leq k\leq N$, let $\beta_k:=s_{i_1}\ldots s_{i_{k-1}}(\alpha_{i_k})$. The sequence $(\beta_1,\beta_2,\ldots,\beta_N)$ is birational.
\item[(3)] If we enumerate all positive roots $\gamma_1,\gamma_2,\ldots,\gamma_N$ in such a way that for $1\leq i,j\leq N$, if $\gamma_i-\gamma_j$ is a positive root, then $i<j$. Following \cite{FFL3}, such an enumeration is called \emph{good}. The sequence $(\gamma_1,\gamma_2,\ldots,\gamma_N)$ is birational.
\end{enumerate}
\end{example}

\subsection{Newton-Okounkov bodies}
Given a birational sequence $S=(\beta_1,\cdots,\beta_N)$, we define a valuation on $\mathbb{C}[G/B]$.

The birational map in the definition of a birational sequence induces a birational map
$$\mathbb{C}^N\to U_{-\beta_1}\times\cdots\times U_{-\beta_N}\dashrightarrow U^-\to G/B,$$
$$(t_1,\ldots,t_N)\mapsto \exp(t_1f_{\beta_1})\cdots\exp(t_Nf_{\beta_N})\cdot B.$$
Passing to the field of rational functions, this gives an isomorphism
$$\mathbb{C}(t_1,\ldots,t_N)\cong \mathbb{C}(G/B).$$
Fixing a monomial order $>$ on $\mathbb{Z}^N$, one defines a valuation $\nu_>:\mathbb{C}[t_1,\ldots,t_N]\setminus\{0\}\to(\mathbb{N}^N,>)$ by sending a polynomial to its lowest exponent with respect to $>$. Extending this valuation by the rule $\nu_>(f/g)=\nu_>(f)-\nu_>(g)$ for $f,g\in\mathbb{C}[t_1,\ldots,t_N]\setminus\{0\}$ gives a valuation
$$\nu_>:\mathbb{C}(G/B)\setminus\{0\}\cong \mathbb{C}(t_1,\ldots,t_N)\setminus\{0\}\to(\mathbb{Z}^N,>).$$

By fixing a highest weight section $s_\lambda\in V(\lambda)^*$, the homogeneous coordinate ring $\mathbb{C}[G/B]$ can be embedded into $\mathbb{C}(G/B)$ by sending a homogeneous function $g\in V(k\lambda)^*$ of degree $k$ to $g/s_\lambda^k$. The theory of Newton-Okounkov bodies associates to such a valuation
\begin{enumerate}
\item a monoid $\Gamma(S,>,\lambda):=\{(k,\mathbf{a})\mid k\in\mathbb{N},\ \mathbf{a}\in\nu_>(V(k\lambda)^*\setminus\{0\})\}\subseteq \mathbb{N}\times\mathbb{Z}^N$;
\item a convex body $\Delta(S,>,\lambda):=\mathrm{cone}(\Gamma(S,>,\lambda))\cap(\{1\}\times \mathbb{R}^N)$, called the Newton-Okounkov body associated to $S$ and $>$.
\end{enumerate}

If the monoid $\Gamma(S,>,\lambda)$ is finitely generated, there exists a flat degeneration of $G/B$ to the toric variety associated to the monoid.

\subsection{Identifications of Newton-Okounkov bodies}
When the birational sequence and the monomial order are well-chosen, one can recover well-known toric degenerations of $G/B$ arising from canonical or global crystal bases, as well as PBW bases.

In practice, since $\mathbb{C}[G/B]$ is graded by the weight lattice with finite-dimensional graded components, the well-ordering property holds automatically; hence it suffices to consider monomial orders.

\begin{definition}
We consider the following four monomial orders on $\mathbb{Z}^N$: for $\ba,\bb\in\mathbb{Z}^N$,
\begin{enumerate}
\item lexicographic order: $\ba>_{\mathrm{lex}}\bb$ if the first non-zero coordinate of $\ba-\bb$ is positive;
\item negative lexicographic order: $\ba>_{\mathrm{neglex}}\bb$ if $\bb>_{\mathrm{lex}}\ba$;
\item degree reverse lexicographic order: $\ba>_{\mathrm{degrevlex}}\bb$ if the total degree $\text{deg } \ba > \text{deg } \bb$ or if they are equal, the first non-zero coordinate of $\ba - \bb$, when read from the right, is negative.
\item weighted degree reverse lexicographic order $>_{\mathrm{wdegrevlex}}$ : similar of $\mathrm{degrevlex}$ but the degree is given for each variable by a weight.
\end{enumerate}
\end{definition}

\begin{theorem}[\cite{FaFoL, FFL3, FN}]\label{thm:essential}
Let $\underline{w}_0=s_{i_1}s_{i_2}\ldots s_{i_N}$ be a reduced decomposition of $w_0$.
\begin{enumerate}
\item Fix the birational sequence $S=(\alpha_{i_1},\alpha_{i_2},\ldots,\alpha_{i_N})$.
\begin{enumerate}
\item[(a)] The associated Newton-Okounkov body $\Delta(S,>_{\mathrm{neglex}},\lambda)$ is the Littelmann-Berenstein-Zelevinsky polytope.
\item[(b)] The Newton-Okounkov body $\Delta(S,>_{\mathrm{degrevlex}},\lambda)$ is the Nakashima-Zelevinsky polytope.
\end{enumerate}
\item Fix the birational sequence $S=(\beta_1,\beta_2,\ldots,\beta_N)$ as in Example \ref{Ex:Bir} (2).
The Newton-Okounkov body $\Delta(S,>_{\mathrm{wdegrevlex}},\lambda)$, where the weight of $\beta$ is given by the height of $\beta$, is the Lusztig polytope parametrizing canonical bases.
\item Let $G=\mathrm{SL}_{n+1}$ or $\mathrm{Sp}_{2n}$, $S$ be a good sequence as in Example \ref{Ex:Bir} (3).
The Newton-Okounkov body $\Delta(S,>_{\mathrm{neglex}},\lambda)$ is the Feigin-Fourier-Littelmann-Vinberg polytope $\mathrm{FFLV}(\lambda)$.
\end{enumerate}
\end{theorem}

As a summary, birational sequence is a tool to unify well-known toric degenerations of flag varieties.
The Conjecture~\ref{cnj:main}, that the special fibre of our degeneration is in fact a toric variety, is translated into checking whether the monoid is finitely generated.
All polytopes in Theorem~\ref{thm:essential} are rational, hence the monoids are finitely generated.

\subsection{Essential monomials}
Computing points in a Newton-Okounkov body is a difficult task. 
As an advantage of the approach via birational sequences, rational points can be calculated using representation theory of Lie algebras.

We fix a birational sequence $S=(\beta_1,\beta_2,\ldots,\beta_N)$ and a monomial order on $\mathbb{N}^N$. For $\mathbf{k}=(k_1,\ldots,k_N)\in\mathbb{N}^N$, we denote $f^{\mathbf{k}}:=f_{\beta_1}^{k_1}\ldots f_{\beta_N}^{k_N}$. These data allow us to define a filtration on the universal enveloping algebra $U(\mathfrak{n}^-)$ as follows: for $\bm\in\mathbb{N}^N$,
$$U(\mathfrak{n}^-)_{\leq \bm}:=\mathrm{span}\{f^{\mathbf{k}}\mid \mathbf{k}\leq\bm\}.$$
We set similarly $U(\mathfrak{n}^-)_{< \bm}:=\mathrm{span}\{f^{\mathbf{k}}\mid \mathbf{k}<\bm\}$.
This filtration on $U(\mathfrak{n}^-)$ induces a filtration on the module $V(\lambda)$ by assigning
$$V(\lambda)_{\leq \bm}:= U(\mathfrak{n}^-)_{\leq\bm}\cdot \bv_\lambda,$$
where $\bv_\lambda$ is a highest weight vector of $V(\lambda)$. We denote $V(\lambda)_{<\bm}:= U(\mathfrak{n}^-)_{<\bm}\cdot \bv_\lambda$.

\begin{definition}
A tuple $\bm\in\mathbb{N}^N$ is called an \emph{essential exponent} with respect to the birational sequence $S$, the monomial order $>$, and the weight $\lambda$, if
$$\dim \left(V(\lambda)_{\leq\bm}/V(\lambda)_{<\bm}\right)=1.$$
The set of all essential exponents will be denoted by $\mathrm{es}(S,>,\lambda)$. For $\bm\in\mathrm{es}(S,>,\lambda)$, the monomial $f^\bm$ is called an \emph{essential monomial}.
\end{definition}

For two dominant integral weights $\lambda,\mu$, a big part of $\mathrm{es}(S, > , \lambda + \mu)$ can be computed from the following property where $+$ stands for the Minkowski sum of sets:
\begin{equation}\label{eq:minkowski}
    \mathrm{es}(S, > , \lambda) + \mathrm{es}(S, > ,  \mu) \subseteq \mathrm{es}(S, > , \lambda + \mu).
\end{equation}

It is shown in \cite{FaFoL} that the set $\{f^{\mathbf{m}}\cdot\bv_\lambda\mid \mathbf{m}\in\mathrm{es}(S,>,\lambda)\}$ forms a vector space basis of $V(\lambda)$, called the essential basis of $V(\lambda)$ with respect to $S$ and $>$. Essential bases are monomial bases for irreducible representations.

\begin{theorem}[\cite{FaFoL}]
The following holds:
$$\bigcup_{k\in\mathbb{N}}\{k\}\times \mathrm{es}(S,>,k\lambda)=\Gamma(S,>,\lambda).$$
\end{theorem}

This theorem is the theoretic cornerstone of our algorithm for computing essential exponents. We illustrate the theorem in the following example.

\begin{example}
Let $G=\mathrm{SL}_3$ and $\g=\mathfrak{sl}_3$ be the Lie algebra of traceless $3\times 3$-matrices with Lie bracket the commutator of matrices. Fix $\mathfrak{n}^-$ to be the subalgebra of strictly lower triangular matrices in $\g$ and $\mathfrak{h}$ to be the diagonal matrices in $\g$. The simple roots are denoted by $\alpha_1,\alpha_2$, and the positive roots are $\alpha_1,\beta:=\alpha_1+\alpha_2,\alpha_2$. The simple reflections $s_1$ and $s_2$  with respect to the simple roots generate the Weyl group of $\g$, which is isomorphic to the symmetric group $\mathfrak{S}_3$. We fix a reduced decomposition of the longest element $\underline{w}_0=s_1s_2s_1$.

We consider the adjoint representation of $\g$ on itself: it is an irreducible representation of highest weight $\beta$. Note that if we denote the fundamental weights by $\varpi_1$ and $\varpi_2$, then $\beta=\varpi_1+\varpi_2$.

The weight space $V(\beta)_{-\beta}$ in $V(\beta)$ of weight $-\beta$ is one-dimensional. We work out in detail the essential monomial which gives a basis of $V(\beta)_{-\beta}$ for different birational sequences and different monomial orders.

\begin{enumerate}
\item $S=(\alpha_1,\beta,\alpha_2)$ is a birational sequence and we fix the monomial order $>_{\mathrm{neglex}}$. All exponents $\mathbf{m}=(m_1,m_2,m_3)$ such that
$f_{\alpha_1}^{m_1}f_{\beta}^{m_2}f_{\alpha_2}^{m_3}\cdot \bv_\beta$ is non-zero in $V(\beta)_{-\beta}$ are: $(1,1,1)$ and $(0,2,0)$. With respect to $>_{\mathrm{neglex}}$
we would choose $f_{\alpha_1}f_{\beta}f_{\alpha_2}$ as an essential monomial, and then $f_{\alpha_1}f_{\beta}f_{\alpha_2}\cdot \bv_\beta$ is a basis of $V(\beta)_{-\beta}$.
\item Choosing $S=(\beta,\alpha_1,\alpha_2)$ as the birational sequence and the monomial order $>_{\mathrm{neglex}}$.
The basis of $V(\beta)_{-\beta}$ is given by the essential monomial $f_\beta^2\cdot \bv_\beta$ of exponent $(2,0,0)$.
\item Choosing $S=(\alpha_1,\alpha_2,\alpha_1)$ as the birational sequence. There exists only one tuple $\bm=(1,2,1)$ such that $f^\bm\cdot\bv_\beta\in V(\beta)_{-\beta}$ is non-zero.
The essential monomial giving a basis of $V(\beta)_{-\beta}$ is $f_{\alpha_1}f_{\alpha_2}^2f_{\alpha_1}\cdot \bv_\beta$ which is independent of the choice of the monomial order.
\end{enumerate}
\end{example}

\section{Algorithm and first example}
The computing of essential monomials is implemented in OSCAR with
\begin{algorithm}
\caption{Compute monomial basis}
\label{compute_monomials}
\SetAlgoLined
\DontPrintSemicolon
\SetKwInOut{Input}{Input}\SetKwInOut{Output}{Output}
\Input{Lie algebra, birational sequence, order, highest weight $\lambda$, essential monomials for fundamentals, read off essentially from \GAP\!\!.}
\Output{Essential monomials for $\mathrm{es}(\lambda)$.}
\BlankLine
\For{all Minkowski sums $\lambda = \mu_1 + \mu_2$}{
    Compute $\mathrm{es}(S,>,\mu_1)$ and $\mathrm{es}(S,>,\mu_2)$ recursively.\;
    Compute the Minkowski sum $\mathrm{es}(S,>,\mu_1) + \mathrm{es}(S,>,\mu_2)$.\;
}
Union all of these Minkowski sums $\mathrm{es}(S,>,\mu_1) + \mathrm{es}(S,>,\mu_2)$.\;
\If{the total size is smaller than the expected dimension (using Weyl's dim. formula)}{
    Identify weight space where monomials are missing (using the Weyl polytope, \polymake and \GAP\!\!).\;
    Find integer solutions of linear equation system to find missing monomials.\;
}
\end{algorithm}

The use of equation (\ref{eq:minkowski}) in the algorithm reduces the run time as it
helps to avoid solving systems of linear equations if possible but computes Minkowski sums of lattice points only.

Note that the validation of the input is not checked in the program, the program does not
check if the provided sequence is in fact birational.
There will be a message, that the sequence is probably not a birational if
(after a reasonable time) no such generating set can be computed.

Having the goal in mind of providing evidence for Conjecture~\ref{cnj:main} that the essential monoid is finitely generated, there is an additional return value that gives actually a
set of highest weights that are sufficient for generating the essential monomials from
smaller highest weights (component-wise when written in terms of fundamental weights). For given $\lambda$,
the program gives a finite set of dominant integral weights $\{ \mu_i \}$
with $\mu_i$ satisfying $\lambda = \sum a_i \mu_i$ where $a_i\in\mathbb{N}$ and
$$\mathrm{es}(S, > , \lambda) = \sum a_i \mathrm{es}(S,>,\mu_i).$$
We store a set of minimal generators with respect to the graded lexicographic ordering on the coefficients of the highest weights (along with smaller, redundant weights that were computed in the process in step \textbf{2}), that is providing a set of weights whose height is minimal and $\mathrm{es}(S, > , \lambda)$ is in the Minkowski sum of the essential monomials.

The algorithm is implemented in \OSCAR for finite-dimensional complex simple Lie algebras of arbitrary type and rank.

The main function of this \OSCAR package is
\begin{center}
    \texttt{basis\_lie\_highest\_weight(type, rank, highest\_weight, birational\_sequence; monomial\_ordering)}
\end{center}

\medskip
We explain the arguments of the function:

\begin{itemize}
    \item \texttt{type} refers to the type of the Lie algebra as a \julia symbol, for example \texttt{:A} or \texttt{:E}.

    \item \texttt{rank} refers to the rank of the Lie algebra, for example for $\lie{sl}_5(\mathbb{C})$ the rank is $4$.

    \item \texttt{highest\_weight} refers to the highest weight of the module, the input has to be given as a coefficient vector in terms of the fundamental weights, for example for $\lie{sl}_5(\mathbb{C})$ and the adjoint representation, the input would be \texttt{[1,0,0,1]}.

    \item \texttt{birational\_sequence} refers to a birational sequence. The input of the roots is by their index in the output of
    \begin{center}\texttt{basis\_lie\_highest\_weight\_operators(type, rank)},\end{center}
    where the roots are written as linear combinations of the simple roots. For example
\begin{minted}{jlcon}
julia> basis_lie_highest_weight_operators(:A, 4)
10-element Vector{Tuple{Int64, Vector{QQFieldElem}}}:
 (1, [1, 0, 0, 0])
 (2, [0, 1, 0, 0])
 (3, [0, 0, 1, 0])
 (4, [0, 0, 0, 1])
 (5, [1, 1, 0, 0])
 (6, [0, 1, 1, 0])
 (7, [0, 0, 1, 1])
 (8, [1, 1, 1, 0])
 (9, [0, 1, 1, 1])
 (10, [1, 1, 1, 1])
\end{minted}
    The birational sequence 
    \begin{gather}\label{eq:operators-example}
       (\alpha_1, \alpha_2, \alpha_3, \alpha_4, \alpha_1, \alpha_1 + \alpha_2, \alpha_1 + \alpha_2 + \alpha_3, \alpha_2, \alpha_2 + \alpha_3, \alpha_3)
    \end{gather} 
    would be translated to \texttt{[1,2,3,4,1,5,8,2,6,3]}.
    
    Alternatively, one can provide a list of coefficient vectors with respect to simple roots. The previous birational sequence would be translated to \texttt{[[1,0,0,0], [0,1,0,0], [0,0,1,0], [0,0,0,1], [1,0,0,0], [1,1,0,0], [1,1,1,0], [0,1,0,0], [0,1,1,0], [0,0,1,0]]}.
    
    As already mentioned, it is not checked by the function that the input is indeed a birational sequence.

    \item \texttt{monomial\_ordering} refers to a monomial ordering, the default is \texttt{degrevlex}. Further monomial orderings can be found in the \OSCAR documentation\footnote{\url{https://docs.oscar-system.org/v1/CommutativeAlgebra/GroebnerBases/orderings/}}.
\end{itemize}

\begin{example}
We use the function for a Lie algebra of type $A$, of rank $4$, for the module of dominant integral weight $2 \varpi_1 + \varpi_2 + 2\varpi_3 + \varpi_4$, the birational sequence from (\ref{eq:operators-example}) and the \textit{degrevlex} ordering.
\begin{minted}{jlcon}
julia> b = basis_lie_highest_weight(:A, 4, [2,1,2,1], [1,2,3,4,1,5,8,2,6,3]; monomial_ordering=:degrevlex)
Monomial basis of a highest weight module
  of highest weight [2, 1, 2, 1]
  of dimension 8750
  with monomial ordering degrevlex([x1, x2, x3, x4, x5, x6, x7, x8, x9, x10])
over Lie algebra of type A4
  where the used birational sequence consists of the following roots (given as coefficients w.r.t. alpha_i):
    [1, 0, 0, 0]
    [0, 1, 0, 0]
    [0, 0, 1, 0]
    [0, 0, 0, 1]
    [1, 0, 0, 0]
    [1, 1, 0, 0]
    [1, 1, 1, 0]
    [0, 1, 0, 0]
    [0, 1, 1, 0]
    [0, 0, 1, 0]
  and the basis was generated by Minkowski sums of the bases of the following highest weight modules:
    [1, 0, 0, 0]
    [0, 1, 0, 0]
    [0, 0, 1, 0]
    [0, 0, 0, 1]
\end{minted}

\bigskip

The output is the monomial basis \texttt{monomials(b)}, even more certain properties of this basis are also computed such as its length. We can read off from
\begin{minted}{jlcon}
  and the basis was generated by Minkowski sums of the bases of the following highest weight modules:
    [1, 0, 0, 0]
    [0, 1, 0, 0]
    [0, 0, 1, 0]
    [0, 0, 0, 1]
\end{minted}
\noindent
that $\mathrm{es}(S, >, 2\varpi_1 + \varpi_2 + 2\varpi_3 + \varpi_4)$ is a Minkowski sum of $\mathrm{es}(S, >, \varpi_1)$, $\mathrm{es}(S, >, \varpi_2)$, $\mathrm{es}(S, >, \varpi_3)$ and $\mathrm{es}(S, >, \varpi_4)$. 
In fact:
$$\mathrm{es}(S,>,2\varpi_1 + \varpi_2 + 2\varpi_3 + \varpi_4) = 2\mathrm{es}(S,>, \varpi_1) + \mathrm{es}(S, >, \varpi_2) + 2\mathrm{es}(S, >, \varpi_3) + \mathrm{es}(S, >, \varpi_4).
$$
\end{example}

\section{More examples}
There are shorthands implemented for important classes of birational sequences and monomial orderings in Theorem~\ref{thm:essential}.

\begin{example}
    Let $\lie g$ be a finite-dimensional complex simple Lie algebra, $S$ be the set of positive roots, ordered by their heights in a descending way and $\lambda$ be a dominant integral weight. Then the essential monomials
    with respect to the monomial order $>_{\text{degrevlex}}$ can be computed by
\begin{minted}{jlcon}
julia> basis_lie_highest_weight_ffl(:A, 3, [1,1,1])
Monomial basis of a highest weight module
  of highest weight [1, 1, 1]
  of dimension 64
  with monomial ordering degrevlex([x1, x2, x3, x4, x5, x6])
over Lie algebra of type A3
  where the used birational sequence consists of the following roots (given as coefficients w.r.t. alpha_i):
    [1, 1, 1]
    [0, 1, 1]
    [1, 1, 0]
    [0, 0, 1]
    [0, 1, 0]
    [1, 0, 0]
  and the basis was generated by Minkowski sums of the bases of the following highest weight modules:
    [1, 0, 0]
    [0, 1, 0]
    [0, 0, 1]
\end{minted}
    For $\mathfrak{g} = \lie{sl}_{n+1}(\CC)$ or $\lie{sp}_{2n}(\CC)$, the bases are known explicitly by \cite{FFL1, FFL2} and the Newton-Okounkov body $\Delta(S,>_{\mathrm{degrevlex}},\lambda)$ is the Feigin-Fourier-Littelmann-Vinberg polytope
    $\mathrm{FFLV}(\lambda)$.
    Again, \polymake can be used to compute the actual FFLV polytope.
    Note that explicit defining inequalities are so far only known for the FFLV polytopes for Lie algebras of types $A, C, G$ (for type $G$ see \cite{Gor15}). In other types, and even for some Lie superalgebras, variants of these FFLV bases are known (\cite{BK19, BD15, Mak19, Gor19, FK21}).
\end{example}

\begin{example}
    Let $\lie g$ be a finite-dimensional complex simple Lie algebra and 
    let $\underline{w}_0=s_{i_1}s_{i_2}\ldots s_{i_N}$ be a
    reduced decomposition of $w_0$, the longest element in the Weyl group of $\mathfrak{g}$.
    We consider the birational sequence $(\alpha_{i_1},\alpha_{i_2},\ldots,\alpha_{i_N})$ 
    and fix as the monomial order $>_{\mathrm{neglex}}$.
    Then the monoid $\Gamma$ is finitely generated and the associated 
    Newton-Okounkov body $\Delta(S,>_{\mathrm{neglex}},\lambda)$
    is the Littelmann-Berenstein-Zelevinsky polytope, also known as the string polytope.
    This very important class of examples can be computed by the shorthand
\begin{minted}{jlcon}
julia> basis_lie_highest_weight_string(:B, 3, [1,1,1], [3,2,3,2,1,2,3,2,1])
Monomial basis of a highest weight module
  of highest weight [1, 1, 1]
  of dimension 512
  with monomial ordering neglex([x1, x2, x3, x4, x5, x6, x7, x8, x9])
over Lie algebra of type B3
  where the used birational sequence consists of the following roots (given as coefficients w.r.t. alpha_i):
    [0, 0, 1]
    ...
    [1, 0, 0]
  and the basis was generated by Minkowski sums of the bases of the following highest weight modules:
    [1, 0, 0]
    [0, 1, 0]
    [0, 0, 1]

\end{minted}
\end{example}

\begin{example}
    Closely related to the string polytope are Lusztig polytopes. 
    We fix the same notation as in the previous example. 
    Then $S = (\beta_1, \ldots, \beta_N)$ with $\beta_k = s_{i_1} \cdots s_{i_{k-1}} (\alpha_{i_k})$ is a birational sequence (Example \ref{Ex:Bir} (2)).
    The weight of $\beta$ is given by the sum of the coefficients in the linear combination of simple roots.
    The Newton-Okounkov body $\Delta(S,>_{\mathrm{wdegevlex}},\lambda)$ is 
    the Lusztig polytope parametrizing the canonical basis of $V(\lambda)$.
\begin{minted}{jlcon}
julia> basis_lie_highest_weight_lusztig(:D, 4, [1,1,1,1], [4,3,2,4,3,2,1,2,4,3,2,1])
Monomial basis of a highest weight module
  of highest weight [1, 1, 1, 1]
  of dimension 4096
  with monomial ordering wdegrevlex([x1, x2, x3, x4, x5, x6, x7, x8, x9, x10, x11, x12], [1, 1, 3, 2, 2, 1, 5, 4, 3, 3, 2, 1])
over Lie algebra of type D4
  where the used birational sequence consists of the following roots (given as coefficients w.r.t. alpha_i):
    [0, 0, 0, 1]
    ...
    [1, 0, 0, 0]
  and the basis was generated by Minkowski sums of the bases of the following highest weight modules:
    [1, 0, 0, 0]
    [0, 1, 0, 0]
    [0, 0, 1, 0]
    [0, 0, 0, 1]
    [0, 0, 1, 1]
\end{minted}
\end{example}

\begin{example}
    A class of examples having very interesting properties is the following.
    We fix again the same notation as in the previous example and 
    the birational sequence $S=(\alpha_{i_1},\alpha_{i_2},\ldots,\alpha_{i_N})$.
    The Newton-Okounkov body $\Delta(S,>_{\mathrm{degrevlex}},\lambda)$ is the 
    Nakashima-Zelevinsky polytope.
    The difference to the string polytope is the chosen order, instead of negative 
    lexicographic one chooses the degree reverse lexicographic order.
\begin{minted}{jlcon}
julia> basis_lie_highest_weight_nz(:C, 3, [1,1,1], [3,2,3,2,1,2,3,2,1])
Monomial basis of a highest weight module
  of highest weight [1, 1, 1]
  of dimension 512
  with monomial ordering degrevlex([x1, x2, x3, x4, x5, x6, x7, x8, x9])
over Lie algebra of type C3
  where the used birational sequence consists of the following roots (given as coefficients w.r.t. alpha_i):
    [0, 0, 1]
    ...
    [1, 0, 0]
  and the basis was generated by Minkowski sums of the bases of the following highest weight modules:
    [1, 0, 0]
    [0, 1, 0]
    [0, 0, 1]
\end{minted}
\end{example}

The Gelfand-Tsetlin polytope for $\lie g = \mathfrak{sl}_{n+1}$ is (up to a linear translation) an instance of a Nakashima-Zelevinsky polytope,
namely for the reduced expression $(s_{1}s_2 \ldots s_{n}) (s_1 s_2 \ldots s_{n-1}) \ldots (s_1s_2)s_1.$

\section{Generators of the monoid for fixed Kodaira embedding}
Attempting to compute $\Gamma(S, >, \lambda)$, and especially assessing its finite generation property, naturally requires computing $\mathrm{es}(S, >, n \lambda)$ for a fixed $\lambda$ and all $n$. Conjecture~\ref{cnj:main} translates to: does there exist an $N > 0$ such that for all $k \geq 0$, $\mathrm{es}(S, >, (N + k) \lambda)$ can be expressed as a Minkowski sum of the sets $\mathrm{es}(S, >, n \lambda)$ for $n < N$?

To address this, we have developed a function in \OSCAR that calculates, for a given $N$, the generators of the truncated monoid represented as: 
$$
\Gamma(S,>,\lambda) / \Gamma_{N+1} \text{ where } \Gamma_{N+1} = \{0 \} \cup \bigcup_{k > N}\{k\}\times \mathrm{es}(S,>,k\lambda).
$$
The main \OSCAR function in this context is
\begin{center}
    \texttt{basis\_coordinate\_ring\_kodaira(type, rank, highest\_weight, degree, birational\_sequence; monomial\_ordering)}
\end{center}\textbf{}
The inputs for \texttt{type}, \texttt{rank}, \texttt{highest\_weight}, \texttt{birational\_sequence}, and \texttt{monomial\_ordering} are the same as in \texttt{basis\_lie\_highest\_weight} and
\begin{itemize}
    \item \texttt{degree} refers to the truncation, e.g., for $\Gamma(S,>,\lambda) / \Gamma_{7}$ the input would be $6$.
\end{itemize}

The output of the function is a list of pairs, where the $k$-th entry contains a monomial basis of $V(k\lambda)$ and the set
\[
    \mathrm{es}(S, >, k \lambda) \setminus \left( \bigcup_{\ell = 1}^{k-1} \mathrm{es}(S, >, \ell \lambda) + \mathrm{es}(S, >, (k - \ell) \lambda)\right),
\]
that is, the elements of $\mathrm{es}(S, >, k \lambda)$ that are not contained in the Minkowski sum of smaller degrees.

\begin{example}
We consider the Lie algebra of type $G$ and rank $2$, $\lambda = \varpi_1$ ($\alpha_1$ being the short simple root), $ S = ( \alpha_1, \alpha_2, \alpha_1 +  \alpha_2, 2 \alpha_1 + \alpha_2, 3 \alpha_1 + \alpha_2, 3 \alpha_1 + 2 \alpha_2)$, the monomial ordering is \textit{invlex} ($x^{\alpha} > x^{\beta} \Leftrightarrow \exists\, 1 \leq i \leq n : \alpha_n = \beta_n, \ldots, \alpha_{i+1} = \beta_{i+1}, \alpha_i > \beta_i$) and the degree is set to $6$. So the function computes the generators of $$\Gamma(S,>,\varpi_1) / \Gamma_{7}.$$
\begin{minted}{jlcon}
julia> bases = basis_coordinate_ring_kodaira(:G, 2, [1,0], 6; monomial_ordering = :invlex)
6-element Vector{Tuple{MonomialBasis, Vector{ZZMPolyRingElem}}}:
 (Monomial basis of a highest weight module with highest weight [1, 0] over Lie algebra of type G2, [1, x1, x3, x1*x3, x1^2*x3, x3*x4, x1*x3*x4])
 (Monomial basis of a highest weight module with highest weight [2, 0] over Lie algebra of type G2, [x4, x1*x4, x4^2, x3*x4^2, x1*x3*x4^2])
 (Monomial basis of a highest weight module with highest weight [3, 0] over Lie algebra of type G2, [x1^2*x4^2, x4^3, x1*x4^3, x4^4, x1*x4^4, x3*x4^4, x5, x2*x5, x1*x2*x5, x1^2*x2*x5, x3^2*x5, x1*x3^2*x5, x3^3*x5, x1*x3^3*x5])
 (Monomial basis of a highest weight module with highest weight [4, 0] over Lie algebra of type G2, [x4^5, x1*x4^5, x4^6, x3^2*x4*x5, x1*x3^2*x4*x5, x3^2*x4^2*x5, x3^3*x4^2*x5])
 (Monomial basis of a highest weight module with highest weight [5, 0] over Lie algebra of type G2, [x1^2*x4^6, x4^7, x1*x4^7, x2*x4^3*x5, x1*x2*x4^3*x5, x2*x3*x4^3*x5, x1*x2*x3*x4^3*x5, x1^2*x2*x3*x4^3*x5, x2*x3^2*x4^3*x5, x1*x2*x3^2*x4^3*x5, x1^2*x2*x3^2*x4^3*x5, x2*x4^4*x5])
 (Monomial basis of a highest weight module with highest weight [6, 0] over Lie algebra of type G2, [x4^9, x1*x3*x4^4*x5, x2*x4^5*x5, x3*x4^5*x5, x3^2*x4^5*x5, x2*x3^2*x4^5*x5, x1*x2*x3^2*x4^5*x5, x3^4*x4*x5^2])

julia> [length(basis[2]) for basis in bases]
6-element Vector{Int64}:
  7
  5
 14
  7
 12
  8
\end{minted}
The numbers $7, 5, 14, 7, 12, 8$ are referring to the number of generators added in each degree.
\end{example}

Similar to \texttt{basis\_lie\_highest\_weight\_*}, there are again shorthands for some classes of birational sequences and monomial orderings.


\section{Experiments on  generators of the monoid}
One of the primary motivations for this implementation is to provide evidence supporting Conjecture~\ref{cnj:main}, which posits that the monoid of essential monomials $\Gamma(S,>, \lambda)$ is finitely generated for all birational sequences (\cite{FaFoL}). According to \cite{BZ}, this holds true for all string parametrizations (or equivalently, for all Lusztig parametrizations). However, they do not offer a set of generators or even an upper bound for the degree of elements in a minimal set of generators. In cases where the conjecture is proven true, an obvious question arises: how can one find a set of generators?
In our implementation, we calculated the minimal building blocks for $\mathrm{es}(S, >, 2\rho)$ using the string parametrization outlined in Theorem\ref{thm:essential} for $\lie{sl}_6(\CC)$ and the highest weight $2\rho$, which represents the sum of all positive roots. Specifically, for every reduced expression of the longest element in $S_6$, we determined the essential monomials and, importantly, the necessary elements in the monoid to express $\mathrm{es}(S, >, 2\rho)$ as a Minkowski sum. To eliminate redundant cases, we utilized \OSCAR to generate a list of reduced expressions, accounting for the exchange of successive orthogonal reflections, and then computed the generators for the remaining 908 cases.
We organized the reduced expressions according to the necessity of generators, and here is the output, correlating each set of generators with the frequency of its occurrence. Evidently, regardless of the reduced expression, all fundamental weights appear in the list of generators. To simplify, we refer to these collectively as ``fundamentals''.

\begin{minted}{jlcon}

[fundamentals] => 320
[fundamentals, [1, 0, 1, 0, 0]] => 70
[fundamentals, [0, 0, 1, 0, 1]] => 70
[fundamentals, [0, 1, 0, 1, 0]] => 54
[fundamentals, [1, 0, 0, 1, 0]] => 28
[fundamentals, [0, 1, 0, 0, 1]] => 28
[fundamentals, [1, 0, 0, 1, 0], [0, 1, 0, 1, 0]] => 46
[fundamentals, [0, 1, 0, 1, 0], [0, 1, 0, 0, 1]] => 46
[fundamentals, [1, 0, 1, 0, 0], [1, 0, 0, 1, 0]] => 30
[fundamentals, [0, 1, 0, 0, 1], [0, 0, 1, 0, 1]] => 30
[fundamentals, [1, 0, 1, 0, 0], [0, 1, 0, 0, 1]] => 8
[fundamentals, [1, 0, 0, 1, 0], [0, 0, 1, 0, 1]] => 8
[fundamentals, [1, 0, 1, 0, 0], [1, 0, 0, 1, 0], [0, 1, 0, 1, 0]] => 70
[fundamentals, [0, 1, 0, 1, 0], [0, 1, 0, 0, 1], [0, 0, 1, 0, 1]] => 70
[fundamentals, [1, 0, 0, 1, 0], [0, 1, 0, 1, 0], [0, 1, 0, 0, 1]] => 10
[fundamentals, [1, 0, 1, 0, 0], [1, 0, 0, 1, 0], [0, 1, 0, 1, 0], [0, 1, 0, 0, 1]] => 10
[fundamentals, [1, 0, 0, 1, 0], [0, 1, 0, 1, 0], [0, 1, 0, 0, 1], [0, 0, 1, 0, 1]] => 10
\end{minted}

Even in the ''worst-case'' scenario, no generator exceeds a height of $2$ for $\mathfrak{sl}_6(\mathbb{C})$ and any arbitrary reduced expression. This observation might even be considered as evidence supporting the bold conjecture that the (global) monoid $\Gamma(S, >) := \sum_{\lambda} \Gamma(S, >, \lambda)$ for string polytopes is not only finitely generated but also has an upper bound on the number of generators.This global monoid parametrizes a basis of the multi-homogeneous coordinate ring. Specifically, the monoid $\Gamma(S, >)$ could be generated by all essential monomials of $V(\mu)$ for $\mu$ such that $\mu \leq 2\rho$.\\

\section{Comparison of the runtime against \GAP}
\GAP provides a function to compute a monomial basis of a highest weight module for simple, complex Lie algebras of any type, as demonstrated in the following code segment: 
\begin{minted}{jlcon}
L:= SimpleLieAlgebra("A", 4, Rationals);; 
V:= HighestWeightModule(L, [1,3,2,1]);; 
\end{minted} 
This function uses the \GAP ordering of positive roots, ranging from simple to highest with respect to the weight, as a birational sequence, with the monomial ordering being \textit{deglex}. Currently, \GAP does not support alterations to either the birational sequence or the monomial ordering. In fact, the basis is constructed simultaneously with the module. For modules of smaller dimensions, the runtime of the two algorithms is comparable. The height of the highest weight is determined by the sum of the coefficients when expressed in terms of fundamental modules. Specifically, if the height of $\lambda$ is $N$, then $V(\lambda)$ represents the Cartan component in an $N$-times tensor product of fundamental representations of the Lie algebra. Our algorithm prefers to avoid computations with matrices when feasible, instead utilizing Minkowski sums of ``smalle'' weights. If the conjecture about the monoid's generators is true, then computing monomial bases for large highest weights simplifies to computing bases for relatively small weights and then calculating their Minkowski sum. This implies that our new algorithm is particularly effective when dealing with large highest weight heights.
We used a \textit{Intel(R) Xeon(R) CPU E5-4617 0  2.90GHz} for the run time comparison depicted in \autoref{tbl:runtime_A3}.
\begin{table}
    \centering
    \begin{tabular}{ccrr}\toprule
        Highest weight & Dimension of module & OSCAR & GAP \\ \midrule
        (1,3,2) & 756 & 0.76& 0.23 \\
        (3,3,3) & 4096 & 0.77& 2.64 \\ 
        (3,4,2)& 4320& 0.90 & 1.78\\
        (4,3,5)& 13500& 0.80& 15.10 \\
        (5,6,3)& 34034& 0.91& 48.74 \\
        (5,5,5)& 46656& 1.18& 187.91 \\
        (5,5,6)& 62244& 1.27& 274.85 \\
        (5,6,5)& 67431& 1.31& 296.22 \\
        (6,5,6)& 82810 & 1.40& 882.80 \\
        (6,6,6)& 117649& 2.70& 1763.19
    \end{tabular}
    \caption{Runtime comparison of type $A$, rank $3$, in seconds}
    \label{tbl:runtime_A3}
\end{table}

In the case of type $A$ and rank $4$, $\rho = (1,1,1,1)$ is a generator of the monoid, so parts of its monomial bases have to be added using \polymake and solving systems of linear equations. We show a runtime comparison with \GAP on the same hardware as above in this case in \autoref{tbl:runtime_A4}.
\begin{table}
    \centering
    \begin{tabular}{ccrr}\toprule
        Highest weight & Dimension of module & OSCAR (Jan 24)  & GAP \\ \midrule 
        (1,3,2,1) & 31185 & 3.6 &11.7 \\
        (1,2,2,3)& 50400& 2.4 & 29.9\\
        (2,2,2,2)& 59049& 12.4 & 64.4\\
        (1,3,3,2) & 160160 & 5.2 & 98.8 \\ 
        (2,3,2,2)& 143325&  15.5 &  175.6 \\
        (3,2,2,2)& 110565 & 16.5 &  206.5 \\
        (2,3,3,2)& 332024& 20.9 &  507.0 
    \end{tabular}
    \caption{Runtime comparison of type $A$, rank $4$, in seconds}
    \label{tbl:runtime_A4}
\end{table}

\FloatBarrier 
\printbibliography

\end{document}